\documentclass[12pt,a4paper]{article}
\usepackage{amsmath,amsfonts,amssymb,amscd}

\usepackage[a4paper]{geometry}

\def\id{\operatorname{id}}

 \newtheorem{thm}{Theorem} 
 
 \newtheorem{lem}[thm]{Lemma}
 \newtheorem{prop}[thm]{Proposition}

 \newtheorem{defn}[thm]{Definition}

 \newtheorem{rem}[thm]{Remark}

\newcounter{th}

\newcounter{le}

\newcounter{lem}

\newcounter{de}

\newcounter{ex}

\usepackage{cite}

\begin{document}

\begin{center}
{\bfseries\LARGE	On $n$-ary Lie algebras of type $(r,l)$}\footnote{ Supported by  AFR-grant, University of Luxembourg.
}
\end{center}

\begin{center}
	E.G. Vishnyakova
\end{center}

\begin{abstract} 
These notes are devoted to the multiple generalization of a Lie algebra introduced by A.M.~Vinogradov and M.M.~Vinogradov. We compare definitions of such algebras in the usual and invariant case. Furthermore, we show that there are no simple $n$-ary Lie algebras of type $(n-1,l)$ for $l>0$.
\end{abstract}

\section{Introduction}
This paper is devoted to the study of the multiple generalization of a Lie algebra introduced in \cite{VV1}. 
More precisely in \cite{VV1} A.M.~Vinogradov and M.M.~Vi\-no\-gradov proposed a two-parameter family of $n$-ary algebras that generalized {\it Filippov $n$-algebras} and {\it Lie $n$-algebras}. V.T.~Filippov \cite{Filippov} considered alternating $n$-ary algebras $V$ satisfying
the following Jacobi identity:
\begin{equation}\label{eq_ Jacobi Filippov}
\{ a_1, \ldots ,a_{n-1}, \{ b_1, \ldots, b_n\}\}
=\sum\limits_{i=1}^n
\{ b_1, \ldots, b_{i-1}\{ a_1, \ldots, a_{n-1}, b_i\}, \ldots, b_n\},
\end{equation}
where $a_i,\,b_j\in V$. In other words, the linear maps $\{ a_1, \ldots ,a_{n-1}, - \}:V\to V$ are
derivations of the $n$-ary bracket $\{ b_1, \ldots, b_n \}$.
Another natural $n$-ary generalization of the standard Jacobi identity has the following form:
\begin{equation}\label{eq_ Jacobi SH}
\sum_{(I,J)} (-1)^{(I,J)} \{ \{ a_{i_1},\ldots, a_{i_n}\}, a_{j_1}, \ldots, a_{j_{n-1}}\} = 0,
\end{equation}
where the sum is taken over all ordered unshuffle multi-indexes
$I=(i_1, \ldots,i_n)$ and $J= (j_1, \ldots, j_{n-1})$  such that
$(I,J)$ is a permutation of the multi-index $(1,\ldots, 2n-1)$. We  call alternating $n$-ary algebras satisfying Jacobi identity (\ref{eq_ Jacobi SH}) {\it Lie $n$-algebras}. 

The $n$-ary algebras of type (\ref{eq_ Jacobi Filippov}) appear naturally in  Nambu
mechanics \cite{Nambu} in the context of Nambu-Poisson manifolds, in supersymmetric
gravity theory and in supersymmetric gauge theories, the Bagger-Lambert-Gustavsson Theory, see \cite{n-ary_algebras:_a_review} for details. The $n$-ary algebras of type (\ref{eq_ Jacobi SH}) were considered for instance by P.~Michor and A.M.~Vinogradov in  \cite{MV} and by P.~Hanlon
and M.L.~Wachs \cite{HW}.
The homotopy case was studied in \cite{SS} in the context of the
Schlesinger-Stasheff homotopy algebras and  $L_{\infty}$-algebras. Such algebras are related to
the Batalin-Fradkin-Vilkovisky theory and to string field theory \cite{Lads_Stasheff}.

The paper is structured as follows. In Section $2$ we remind the definition of the Nijenhuis-Richardson bracket on $\bigwedge V^*\otimes V$, where $V$ is a finite dimensional vector space.
In Section $3$ we give a definition Vinogradovs' algebras. In Section $4$ we compare definitions of usual Vinogradovs' algebras and Vinogradovs' algebras with a symmetric non-degenerate invariant form. In Section $4$ we show that there are no simple $n$-ary Lie algebras of type $(n-1,l)$, $l>0$. The case of  $n$-ary Lie algebras of type $(n-1,0)$ was studied in \cite{Ling}.

\section{Nijenhuis-Richardson bracket}

In this section we follow \cite{VV1}. 
 Let $V$ be a finite dimensional vector space over $\mathbb R$ or $\mathbb C$. 
 We put $I^{n}:= (1,\ldots,n)$ and we denote by $I$ and $J$ the ordered multi-indexes $I= (i_{1},\ldots, i_l)$ and $J= (j_{1},\ldots, j_k)$ such that $i_{1}<\ldots< i_l$ and $j_{1}<\ldots < j_k$. Let $|I|:=l$ and $|J|: = k$. Assume that $i_p\ne j_q$ for all $p,q$. Then we have the non-ordered multi-index 
 $$
 (I,J):= (i_{1},\ldots, i_l,j_{1},\ldots, j_k).
 $$
  We denote by $I+J$ the ordering of the multi-index $(I,J)$ and by $(-1)^{(I,J)}$ the parity of the permutation that maps $(I,J)$ into $I+J$. We set
 \begin{align*}
  a_{I^n}:= (a_1\ldots, a_n)\quad a_{I}:= (a_{i_1}\ldots, a_{i_l}), \quad a_{J}:= (a_{j_1}\ldots, a_{j_k}).
    \end{align*}
Let us take $L\in \bigwedge^l V^*\otimes V$ and $K\in \bigwedge^k V^*\otimes V$.
We define:
\begin{equation}\label{eq L[K]}
L[K](a_{I^{l+k-1}} ) = \sum_{
{\scriptstyle I+J=I^{l+k-1}}} (-1)^{(I,J)} L(K(a_{I}), a_{J}).
\end{equation} 
Here we assume that $|I|=k,\, |J|= l-1$.

\begin{defn}\label{de Richardson-Nijenhuis brackets even} The {\it Nijenhuis-Richardson bracket} 
$$
[L,K]^{NR}\in \bigwedge^{l+k-1} V^*\otimes V,
$$
of $L\in \bigwedge^l V^*\otimes V$ and $K\in \bigwedge^k V^*\otimes V$ is defined by the following formula:
$$
[L,K]^{NR}:= (-1)^{(l-1)(k-1)} L[K] - K[L].
$$
\end{defn}

For $L\in \bigwedge^l V^*\otimes V$, we denote by $L_{a_1,\ldots,a_p}$ the following $(l-p)$-linear map:
$$
L_{a_1,\ldots,a_p}(b_1,\ldots, b_{l-p}):= L(a_1,\ldots,a_p,b_1,\ldots, b_{l-p}), 
$$
where $a_i, b_j\in V$.

\section{Vinogradovs' algebras}

In this subsection we recall the definition of $n$-ary algebras introduced by A.M.~Vinogradov and M.M.~Vinogradov in \cite{VV1}. The graded versions of such algebras are studied in \cite{VV2}. These $n$-ary algebras are a generalization of $n$-ary Lie algebras of type (\ref{eq_ Jacobi Filippov}) and of type (\ref{eq_ Jacobi SH}) for even $n$. 
  Let $V$ be a finite dimensional complex or real  vector space. 
\begin{defn}\label{def n-ary algebra}
An {\it $n$-ary algebra structure} $\mu$ on $V$ is an $n$-linear map:
$$
\mu :\underbrace{V\times \cdots \times V}_{n\,\, \text{times}} \to V.
$$
An {\it $n$-ary algebra} is a pair $(V,\mu)$. 
\end{defn}

Sometimes we will use the bracket notation $\{a_1,\ldots, a_n \}$ instead of $\mu (a_1,\ldots, a_n)$. 

\begin{defn}\label{def (skew)symmetric n-ary algebra}  An $n$-ary algebra is called  {\it skew-symmetric} if 
\begin{equation}\label{eq skew-symmetric algebras}
\{a_1,\ldots, a_i, a_{i+1}, \ldots, a_n \} =  -\{a_1,\ldots, a_{i+1}, a_i, \ldots, a_n\}
\end{equation}
for any $a_p\in V$.

\end{defn}

\begin{defn}[\cite{VV1}]\label{def Vinogradov algebras RN} Let $(V,\mu)$ be an $n$-ary skew-symmetric algebra. We say that $(V,\mu)$ is an {\it $n$-ary Lie algebra of type $(r,l)$} if the following holds:
\begin{equation}
\label{eq Vinogradov algebras}
[\mu_{a_1,\ldots, a_r}, \mu_{b_1,\ldots, b_l}]^{NR} = 0,
\end{equation}
where $a_i, b_j\in V$ and $0\leq l \leq r<n$.
\end{defn}

\begin{ex}
Filippov $n$-ary algebras or $n$-ary algebras satisfying (\ref{eq_ Jacobi Filippov}) are exactly $n$-ary Lie algebras of type $(n-1,0)$. For even $n$, the $n$-ary Lie algebras satisfying (\ref{eq_ Jacobi SH}) are exactly $n$-ary algebras of type $(1,0)$. (See \cite{VV1} for details.)
\end{ex}

 \begin{rem}
The theory of Filippov $n$-ary algebras is relatively well-deve\-loped. For instance,
there is a classification of simple real and complex Filippov $n$-ary algebras and an analog of the Levi decomposition \cite{Ling}. W.X.~Ling in \cite{Ling} proved that there exists only one simple finite-dimensional $n$-ary Filippov algebra over an algebraically closed field of characteristic $0$ for any $n>2$.
The simple Filippov $n$-ary superalgebras in the finite and infinite dimensional
case were studied in \cite{Kac}. It was shown there that there are no simple linearly
compact $n$-ary Filippov superalgebras which are not $n$-ary Filippov algebras, if $n > 2$. A
 classification of linearly compact $n$-ary Filippov algebras was given in \cite{Kac}.
 \end{rem}

 \begin{rem}
 The $n$-ary algebras of type (\ref{eq_ Jacobi SH}) were studied for instance in \cite{MV} and \cite{HW}. In \cite{Vi}  the $n$-ary algebras of type (\ref{eq_ Jacobi SH}) endowed with a non-degenerate invariant form were considered.  
  \end{rem}

\section{Nijenhuis-Richardson and Poisson brackets on $V$}

Let $V$ be a finite dimensional vector space that is supplied with a non-degenerate symmetric bilinear form $(\,,)$.
Then $\bigwedge V$ possesses a natural
structure of a Poisson superalgebra defined by the following formulas:
\begin{gather*}
[x,y]: = (x,y), \quad x,y\in W;\\
[v,w_1\cdot w_2]: = [v,w_1]\cdot w_2 + (-1)^{vw_1} w_1\cdot [v,w_2],\\
[v,w] = -(-1)^{vw}[w,v],
\end{gather*}
where $v,w, w_i$ are homogeneous elements in $\bigwedge W$. (See for example \cite{KosSter} or \cite{Vi} for details.) Let us take any element $\tilde{\mu} \in \bigwedge^{n+1} V$. Then we can define an n-ary algebra structure on $V$ in the following way:
\begin{equation}\label{eq der product}
\{a_{1}, \ldots, a_{n}\}: = [a_1,[\ldots,[a_n,\tilde{\mu}]\ldots]], \,\,\, a_i\in V.
\end{equation}

\begin{defn}
A skew-symmetric $n$-ary  algebra structure is called {\it invariant with respect to the form} $(\,,)$ if  the following holds:
\begin{equation} \label{eq quadratic}
(b, \{a_{1}, \ldots, a_{n}\}) = -(a_1, \{b,a_2, \ldots, a_{n}\})
\end{equation}
for any $b, a_{i}\in V$.
\end{defn}

The following observation was noticed in \cite{Vi}:

\begin{prop}
\label{Main obseravation graded}  Assume that $V$ is finite dimensional and $(\,,)$ is
non-degene\-rate.  Any  skew-symmetric invariant $n$-ary algebra structures can be obtained by construction $(\ref{eq der product})$. In other words, for any skew-symmetric invariant $n$-ary algebra structure $\mu\in \bigwedge^nV^*\otimes V$, there exists $\tilde{\mu}\in \bigwedge^{n+1}V$ such that the map $\mu$ coincides with the map defined by the formula $(\ref{eq der product})$.
\end{prop}

Assume that $(V,\mu)$ is an $n$-ary Lie algebra of type $(r,l)$ and the algebra structure $\mu$ is invariant with respect to $(\,,)$. It follows from Proposition \ref{Main obseravation graded} that the algebra $(V,\mu)$ has another $n$-ary algebra structure $\tilde{\mu}\in  \bigwedge^{n+1}(V)$. In this section we study the relationship between these two $n$-ary algebra structures.
Since $\mu$ and $\tilde{\mu}$ define the same algebra structure on $V$, we have:
$$
\mu(a_1,\ldots, a_l) = [a_1,\ldots, [a_l,\tilde{\mu}]\ldots].
$$

Let us take two skew-symmetric invariant $n$-ary algebra structures $L\in \bigwedge^l V^*\otimes V$ and $K\in \bigwedge^k V^*\otimes V$.  It follows from (\ref{eq L[K]}) that
\begin{equation}\label{eq_L[K]}
L[K](a_{I^{l+k-1}}) = 
\sum_{I+J=I^{l+k-1}
} (-1)^{(I,J)}[[a_I, \tilde{K}], [a_J, \tilde{L}]],
\end{equation}
where $|I|=k,\, |J|= l-1$ and for simplicity we denote 
$$
[a_{I},\tilde{L}]:= [a_{i_1},\ldots [a_{i_l},\tilde{L}]\ldots].
$$

\begin{prop}\label{prop RN zamenjaem na Poisson} Let us take two skew-symmetric invariant $n$-ary algebra structures $L\in \bigwedge^l V^*\otimes V$ and $K\in \bigwedge^k V^*\otimes V$. Then we have:
$$
[a_{I^{l+k-1}},[\tilde{L},\tilde{K}]] = - (-1)^{(k+1)(l+1)} [L,K]^{NR}(a_{I^{l+k-1}}).
$$
\end{prop}

\noindent{\bf Proof.} We set 
\begin{align*}
a_{I^{p}}:= (a_1,\ldots, a_{p});\quad
[a_{I^{p}},\tilde {\mu}]:= [a_1,\ldots[a_p,\tilde {\mu}]].
\end{align*} 
Furthermore,
 we have:
\begin{multline*}
[a_{I^{l+k-1}},[\tilde{L},\tilde{K}]]
 =
 \sum\limits_{
 I+J=I^{l+k-1}
 } 
 (-1)^{(I,J) + (k+1)(l+1)}
  [[a_I, \tilde{L}], [a_J, \tilde{K}]]  \\
 +\sum\limits_{ 
 I'+J'=I^{l+k-1}}
 (-1)^{(I',J')+(l+1)k}
 [[a_{I'}, \tilde{L}], [a_{J'}, \tilde{K}]]  \\
  =\sum\limits_{
I+J=I^{l+k-1}
} 
(-1)^{(I,J) + (k+1)(l+1)}
 [[a_I, \tilde{L}], [a_J, \tilde{K}]] \\
+\sum\limits_{ 
I'+J'=I^{l+k-1}}
(-1)^{(J',I')}
[[a_{I'}, \tilde{L}], [a_{J'}, \tilde{K}]]
\\
=\sum\limits_{
 I+J=I^{l+k-1}
} 
(-1)^{(I,J) + (k+1)(l+1)}
 [[a_I, \tilde{L}], [a_J, \tilde{K}]]
 \\
+ \sum\limits_{ I'+J'=I^{l+k-1}}
-(-1)^{(J',I')}
[[a_{J'}, \tilde{K}],[a_{I'}, \tilde{L}]],
\end{multline*}
where $|I|= l$, $|J|= k-1$ and $|I'|= l-1$, $|J;|= k$.
By equation (\ref{eq_L[K]}), we get
\begin{multline*}
[a_{I^{l+k-1}},[\tilde{L},\tilde{K}]]
 =
 (-1)^{(k+1)(l+1)} K[L] (a_{I^{l+k-1}}) -  L[K](a_{I^{l+k-1}}) \\
= - (-1)^{(k+1)(l+1)} [L,K]^{NR}(a_{I^{l+k-1}}).
\end{multline*}
 The proof is complete.$\Box$

 Let $L\in \bigwedge^t V^*\otimes V$ be a skew-symmetric invariant $n$-ary algebra structures. Clearly,
 $$
 L_{a_1,\ldots,a_p}= c [a_1,\ldots,[a_p,\tilde{L}]], 
 $$
where $c\ne 0$. Proposition \ref{prop RN zamenjaem na Poisson} implies the following. If a skew-symmetric invariant algebra $(V,\mu)$ is given by
 an identity of type 
 $$
 [\mu_{a_1,\ldots,a_r},\mu_{a_1,\ldots,a_l}]^{NR}=0,
 $$
 then the same algebra can be given by the following relation:
 $$
  [[a_{I^{r}}\tilde{\mu}],[a_{I^{l}}\tilde{\mu}]]=0.
 $$
The result of our study is the following.

\begin{thm}\label{thm RN and Poisson}
Let $V$ be a finite dimensional vector space with a non-degenerate symmetric bilinear form $(\,,)$ and $\mu \in \bigwedge^n V^*\otimes V$ be an invariant algebra structure on $V$. Then $(V,\mu)$ is an $n$-ary Lie algebra of type $(r,l)$ if and only if the following holds:
$$
 [[a_{I^{r}}\tilde{\mu}],[a_{I^{l}}\tilde{\mu}]]=0
$$ 
for all $a_i,b_j\in V$. 
\end{thm}

\section{Algebras of type $(n-1,l)$}

In this section we consider $n$-ary algebras of type $(n-1,l)$, where $0<l\leq n-1$. More precisely, we study  $n$-ary skew-symmetric algebras satisfying (\ref{eq Vinogradov algebras}) for $k=n-1$ and $l>0$. 

\begin{lem}
Jacobi identity (\ref{eq Vinogradov algebras}) of type $(n-1,l)$, where $0\leq l\leq n-1$, is equivalent to 
the following identity:
\begin{equation}\label{eq Jacobi n-1, l}
\{ a_1,\ldots, a_{n-1},\{ b_1,\ldots, b_{n}\} \} = \sum\limits_{i=1}^{n-l} \{b_1,\ldots, \{a_1,\ldots, a_{n-1},  b_{i}\},\ldots, b_{n} \}.
\end{equation}
\end{lem}

\noindent{\bf Proof.}

We will need the following formula, see \cite{VV1}:
$$
[L,K]^{NR}_a=  [L_a,K]^{NR} + (-1)^{l-1} [L,K_a]^{NR},
$$
where $L\in \bigwedge\limits^l V^*\otimes V$ and $K\in \bigwedge\limits^k V^*\otimes V$. Let $(V,\mu)$ is an $n$-ary Lie algebra of type $(n-1,l)$. Let us take $a_i,b_j,c_t\in V$ and denote $L:= \mu_{a_1,\ldots, a_n}$ and $K:= \mu_{b_1,\ldots, b_{l}}$. We have
$$
[L,K]^{NR}_{c_1,\ldots, c_{k}} = \sum_i (-1)^{k-i} [L_{c_i}, K_{c_1,\ldots,\hat{c}_i,\ldots,c_k} ]^{NR} + [L, K_{c_1,\ldots,c_i,\ldots,c_k} ]^{NR}
$$
where $k=n-l$. Further,
\begin{align*}
[L, K_{c_1,\ldots,c_i,\ldots,c_k} ]^{NR} &= \{a_1,\ldots, a_{n-1}, \{b_1,\ldots, b_{l}, c_1,\ldots, c_{k}\}\};\\
[L_{c_i}, K_{c_1,\ldots,\hat{c}_i,\ldots,c_k} ]^{NR}& = - \{b_1,\ldots, b_{l}, c_1,\ldots,\hat{c}_i,\ldots,c_k, \{a_1,\ldots, a_{n-1}\}\}.
\end{align*}
It is clear that $[L,K]^{NR}=0 $ if and only if $[L,K]^{NR}_{c_1,\ldots, c_{k}} = 0$. The last equality is equivalent to:
\begin{align*}
\{a_1,\ldots, a_{n-1}, \{b_1,\ldots, b_{l}, c_1,\ldots, c_{k}\}\} =\\
 \sum_{i=1}^k (-1)^{k-i} \{b_1,\ldots, b_{l}&, c_1,\ldots,\hat{c}_i,\ldots,c_k, \{a_1,\ldots, a_{n-1}\}\},\\
\{a_1,\ldots, a_{n-1}, \{b_1,\ldots, b_{l}, c_1,\ldots, c_{k}\}\} =\\
 \sum_{i=1}^k  \{b_1,\ldots, b_{l}&, c_1,\ldots,\{a_1,\ldots, a_{n-1}\},\ldots,c_k\},\\
\{a_1,\ldots, a_{n-1}, \{c_1,\ldots, c_{k},b_1,\ldots, b_{l}\}\} =\\
 \sum_{i=1}^k  \{c_1,\ldots,&\{a_1,\ldots, a_{n-1}\},\ldots,c_k,b_1,\ldots, b_{l}\}.
\end{align*}
The proof is complete.$|box$

\begin{ex}
For instance, $2$-ary Lie algebras of type $(1,1)$ are exactly the associative skew-symmetric algebras. Using the fact that these algebras are skew-symmetric, it is easy to see that all such algebras have the following property: $\{\{a,b\},c\} =0$. Therefore, this class is not reach. 
\end{ex}

\begin{defn}
 A vector subspace $W\subset V$ is called an {\it ideal} of a skew-symmetric $n$-ary algebra $(V,\mu)$ if $\mu(V, \ldots, V, W) \subset W$. 
\end{defn}

For example $\{0\}$ and $V$ are always ideals of $(V,\mu)$. Such ideals are called {\it trivial}.

\begin{defn}
 An $n$-ary Lie algebra of type $(n-1,l)$ is called {\it simple} if it is not $1$-dimensional and it does not have any non-trivial ideals.
\end{defn}

\begin{thm}\label{prop net prostyh} Let $V$ be a Lie algebra of type  $(n-1,l)$, where $0<l\leq n-1$, over $\mathbb C$. Assume that $V$ does not possess any non-trivial ideals. Then $\dim V= 1$. 
In other words, there are no simple algebras of type  $(n-1,l)$, where $0<l\leq n-1$.

\end{thm}

\noindent{\bf Proof.} Assume that $\dim V>1$ and $V$ does not possess any non trivial ideals. Let us interchange $b_{n-l}$ and $b_{n-l+1}$ and use (\ref{eq Jacobi n-1, l}). We get:
\begin{multline*}
\{ a_1,\ldots, a_{n-1},\{ b_1,\ldots,b_{n-l-1}, b_{n-l+1} , b_{n-l},b_{n-l+2},\ldots, b_{n}\} \} \\ =\sum\limits_{i=1}^{n-l-1} \{b_1,\ldots, \{a_1,\ldots, a_{n-1},  b_{i}\},\ldots, b_{n} \} \\
 +\{b_1,\ldots,b_{n-l-1}, \{a_1,\ldots, a_{n-1},  b_{n-l+1}\},b_{n-l}, b_{n-l+2},\ldots, b_{n} \}.
\end{multline*}
From this equation and (\ref{eq Jacobi n-1, l}), we have:
\begin{multline}\label{eq col from Jacobi n-1, l}
\{b_1,\ldots, b_{n-l-1},  \{a_1,\ldots, a_{n-1},  b_{n-l}\}, b_{n-l+1}\ldots, b_{n} \} \\
 =\{b_1,\ldots, b_{n-l}, \{a_1,\ldots, a_{n-1},  b_{n-l+1}\}, b_{n-l+2},\ldots, b_{n} \}. 
\end{multline}
Denote by $\mathcal L$ the vector subspace in $\operatorname{End}(V)$ generated by $\mu_{a_1,\ldots, a_{n-1}}$, where $a_i\in V$. Let $D=\mu_{a_1,\ldots, a_{n-1}}$ for some $a_i\in V$. We can rewrite (\ref{eq col from Jacobi n-1, l}) in the form:
\begin{multline}\label{eq D(b_n-1)}
\{b_1,\ldots, b_{n-l-1}, D (b_{n-l}), b_{n-l+1},\ldots, b_{n}\} \\
 =\{b_1,\ldots, b_{n-l}, D( b_{n-l+1}), b_{n-l+2},\ldots, b_{n}\}. 
\end{multline}
The $n$-ary algebra $(V,\mu)$ is skew-symmetric and the equation (\ref{eq D(b_n-1)}) holds for any $b_i$. Hence, we can assume that the following equation holds 
$$
\{b_1,\ldots, D (b_{i}), b_{i+1},\ldots, b_{n}\} = \{b_1,\ldots, b_{i}, D( b_{i+1}),\ldots, b_{n}\} 
$$
for any $i$. Applying this formula several times, we get:
$$
\{b_1,\ldots, D (b_{i}), b_{i+1},\ldots, b_{n}\} = \{b_1,\ldots, b_{i+j-1}, D( b_{i+j}),\ldots, b_{n}\} 
$$
for all $i,\,j$. Furthermore, let us take $D_1,D_2\in \mathcal L$. 
 First of all, we see:
\begin{multline*}
\{b_1,\ldots, D_1(D_2 (b_{i})), b_{i+1},\ldots, b_{n}\} =
 \{b_1,\ldots,D_2 (b_{i}), \ldots,  D_1( b_{i+j}),\ldots, b_{n}\}\\ =
  \{b_1,\ldots,b_{i}, \ldots,  D_2 (D_1( b_{i+j})),\ldots, b_{n}\}.
\end{multline*}
We put $B:=D_1\circ D_2 - D_2\circ D_1.
$
Therefore, 
$$
\{b_1,\ldots, B (b_{i}), \ldots,b_{i+j}, \ldots, b_{n}\} = -\{b_1,\ldots, b_{i}, \ldots, B( b_{i+j}),\ldots, b_{n}\}. 
$$
Furthermore,
\begin{multline*}
\{B(b_1), b_2,b_3,\ldots,  b_{i+1},\ldots, b_{n}\} = -\{b_1, B(b_2), b_3,\ldots,  b_{i+1},\ldots, b_{n}\} \\
=\{b_1, b_2,B(b_3),\ldots,  b_{i+1},\ldots, b_{n}\};
\end{multline*}
\begin{equation*}
\{B(b_1), b_2,b_3,\ldots,  b_{i+1},\ldots, b_{n}\}=-
\{b_1, b_2,B(b_3),\ldots,  b_{i+1},\ldots, b_{n}\}.
\end{equation*}
Therefore, 
$$
\{B(b_1), b_2,b_3,\ldots,  b_{n}\}=0
$$
 for all $b_i\in V$. If $B(b_1)\ne 0$, we see that $\langle B(b_1) \rangle$ is a non-trivial $1$-dimensional ideal in $V$. Hence, $B(b_1)=0$ for all $b_1\in V$. Therefore, $D_1\circ D_2 - D_2\circ D_1 =0$ and $\mathcal L$ is a usual abelian Lie algebra. Since, $V$ does not possess any non-trivial ideals, the $\mathcal L$-module $V$ should be irreducible. By Schur's Lemma $\mathcal L \subset \{ a \id\}\subset \operatorname{End}\, V$ and $\dim V=1$.$\Box$

\begin{ex} The classification of simple complex and real Filippov  $n$-ary
algebras or Lie algebras of type $(n-1,0)$ was done in \cite{Ling}: there is one series of complex
Filippov $n$-ary algebras $A_k$, where $k$ is a natural number, and
several real forms for each $A_k$. All these algebras have invariant
forms and in terms of Poisson bracket they are given by the top form of $V$.

\end{ex}

\noindent{\it Elizaveta Vishnyakova}

\noindent {Max Planck Institute for Mathematics Bonn and}

\noindent University of Luxembourg

\noindent{\emph{E-mail address:}
	\verb"VishnyakovaE@googlemail.com"}

\end{document}